\documentclass[a4paper,12pt]{article}
\begin{document}

\title{On structure of topological metagroups.}
\author{Sergey Victor  Ludkowski}
\date{28 October 2021}
\maketitle
\begin{abstract}
In this article topologies on metagroups are studied. They are
related with generalized $C^*$-algebras over ${\bf R}$ or ${\bf C}$.
Homomorphisms and quotient maps on them are investigated. Structure
of topological metagroups is scrutinized. In particular, topologies
on smashed products and smashed twisted wreath products of
metagroups are scrutinized, which are making them topological
metagroups. Moreover, their inverse homomorphism systems are
studied. \footnote{key words and phrases: metagroup; topology;
inverse homomorphism system; generalized $C^*$-algebra; product;
smashed; twisted wreath   
\par S.V.Ludkowski's ORCID: 0000-0002-4733-8151
\par Address for correspondence: S.V. Ludkowski, Dep. Appl. Mathematics, MIREA - Russian Technological
University, av. Vernadsky 78, Moscow 119454, Russia;
\par e-mail: sludkowski@mail.ru }

\end{abstract}

\section{Introduction.}
\par A topological group structure plays very important role in
mathematics and its applications
\cite{archtkachb,hewross,montzip55,pontr}.

\par On the other side, noncommutative analysis is a very important part of
mathematical analysis and it interacts with operator theory,
operator algebras and algebraic analysis. In particular, analysis
over quaternions, octonions and generalized Cayley-Dickson algebras
develops fast in recent years (see
\cite{baez,dickson,frenludkfejms18}-\cite{guertzeb,kansol}-\cite{ludkcvee13,schaeferb}
and references therein). It also plays a huge role in mathematics,
analysis of noncommutative geometry, mathematical physics, quantum
field theory and quantum gravity, partial differential equations
(PDEs), particle physics, operator theory, etc.
\par It appears that a multiplicative law of their
canonical bases is nonassociative and leads to a more general notion
of a metagroup instead of a group
\cite{dickson,ludlmla18,ludaxiomchna}. They were used in
\cite{ludlmla18,luoradaca10,lustsdoadaca11,ludaaca13,ludaxiomchna}
for investigations of partial differential operators and other
unbounded operators over quaternions and octonions, also for
automorphisms, derivations and cohomologies of generalized
$C^*$-algebras over ${\bf R}$ or ${\bf C}$. They certainly have a
lot of specific features in their derivations and (co)homology
theory \cite{ludlmla18,ludaxiomchna}. It was shown in
\cite{ludaaca13} that an analog of the Stone theorem for one
parameter groups of unitary operators for the generalized
$C^*$-algebras over quaternions and octonions becomes more
complicated and multiparameter. The generalized $C^*$-algebras arise
naturally while decompositions of PDEs or systems of PDEs of higher
order into PDEs or their systems of order not higher than two
\cite{frenludkfejms18,guertzeb,ludkrimut14,ludkcvee13,ludcmft12,ludcma13},
that permits to integrate them subsequently or simplify their
analysis. \par In \cite{ludstwpmgax19} different types of products
of metagroups were studied such as smashed products and smashed
twisted wreath products. That permitted also to construct their
abundant families different from groups. But topologies on such
metagroups were not investigated. We recall their definition.

\par {\bf Definition 1.1.}  Let $G$ be a set with a single-valued
binary operation (multiplication)  $G^2\ni (a,b)\mapsto ab \in G$
defined on $G$ satisfying the conditions: \par $\quad (1.1.1)$ for
each $a$ and $b$ in $G$ there is a unique $x\in G$ with $ax=b$  \par
$\quad (1.1.2)$ and a unique $y\in G$ exists satisfying $ya=b$,
which are denoted by $x=a\setminus b=Div_l(a,b)$ and
$y=b/a=Div_r(a,b)$ correspondingly,

\par $\quad (1.1.3)$ there exists a neutral (i.e. unit) element $e_G=e\in G$: \par $~eg=ge=g$
for each $g\in G$.
\par If the set $G$ with the single-valued
multiplication satisfies the conditions $(1.1.1)$ and $(1.1.2)$,
then it is called a quasigroup. If the quasigroup $G$ satisfies also
the condition $(1.1.3)$, then it is called an algebraic loop (or
unital quasigroup, or shortly a loop).
\par The set of all elements $h\in G$
commuting and associating with $G$:
\par $\quad
(1.1.4)$ $Com (G) := \{ a\in G: \forall b\in G, ~ ab=ba \} $,
\par $\quad (1.1.5)$ $N_l(G) := \{a\in G: \forall b\in G, \forall c\in G, ~ (ab)c=a(bc) \}
$,
\par $\quad (1.1.6)$ $N_m(G) := \{a\in G: \forall b\in G, \forall c\in G, ~ (ba)c=b(ac)
\} $,
\par $\quad (1.1.7)$ $N_r(G) := \{a\in G: \forall b\in G, \forall c\in G, ~ (bc)a=b(ca)
\} $,
\par $\quad (1.1.8)$ $N(G) := N_l(G)\cap N_m(G)\cap N_r(G)$;  \par ${\cal C}(G) := Com (G)\cap N(G)$
is called the center ${\cal C}(G)$ of $G$.
\par We call $G$ a metagroup if a set $G$ possesses a single-valued binary operation
and satisfies the conditions $(1.1.1)$-$(1.1.3)$ and
\par $\quad (1.1.9)$ $(ab)c=t(a,b,c)a(bc)$   \\ for each
$a$, $b$ and $c$ in $G$, where $t(a,b,c)=t_G(a,b,c)\in {\cal C}(G)$.

\par Then the metagroup $G$ will be called a central metagroup,
if it satisfies also the condition:
\par $\quad (1.1.10)$  $ab={\sf t}_2(a,b)ba$ \\ for each $a$ and $b$ in $G$, where
${\sf t}_2(a,b)\in {\cal C}(G)$.
\par If $H$ is a submetagroup (or a subloop) of the metagroup $G$
(or the loop $G$) and \par $\quad (1.1.11)$ $gH=Hg$ for each $g\in
G$, then $H$ will be called almost invariant (or algebraically
almost normal). If in addition \par $\quad (1.1.12)$ $(gH)k=g(Hk)$
and $k(gH)=(kg)H$ for each $g$ and $k$ in $G$, \\ then $H$ will be
called an invariant (or algebraically normal) submetagroup (or
subloop respectively).
\par Henceforward it will be used a notation $Inv_l(a)=Div_l(a,e)$
and $Inv_r(a)=Div_r(a,e)$.
\par Elements of a metagroup $G$ will be denoted by small letters,
subsets of $G$ will be denoted by capital letters. If $A$ and $B$
are subsets in $G$, then $A-B$ means the difference of them $A-B=\{
a\in A: ~a \notin B \} $. Henceforward, maps and functions on
metagroups are supposed to be single-valued if something other will
not be specified.
\par If ${\cal T}_G$ is a topology on the metagroup
(or quasigroup, or loop) $G$ such that multiplication, $Div_l$ and
$Div_r$ are (jointly) continuous from $G\times G$ into $G$, then
$(G, {\cal T}_G)$ is called a topological metagroup (or quasigroup,
or loop respectively).
\par Notice also that a loop is a quite different object than a loop
group considered in geometry or mathematical physics. Certainly
loops are more general objects, than metagroups. Note that
metagroups are commonly nonassociative and having many specific
features in comparison with loops and groups. On the other hand, if
a loop $G$ is simple, then a subloop generated by all elements of
the form $((ab)c)/(a(bc))$ for all $a$, $b$, $c$ in $G$ coincides
with $G$ \cite{alb1trams43,bruckb}. Metagroups are intermediate
between groups and quasigroups. Functions on topological metagroups
were studied in \cite{ludaxsftm20}.
\par We recall, that according to Chapter 2 and Sections 4.6, 4.10, 4.13
in \cite{montzip55} and Section 6 in \cite{hewross} the compact
connected $T_0$ topological group $G$ can be presented as the limit
of an inverse homomorphism system (that is, projective limit)
$G=\lim \{ G_j, \pi ^j_k, \Omega \} $ of compact finite-dimensional
Lie groups of manifolds over ${\bf R}$, where $\Omega $ is a
directed set, $\pi ^j_k: G_j\to G_k$ is a continuous homomorphism
for each $j>k$ in $\Omega $, $\pi ^j_j$ is the identity map, $\pi
^j_j(g_j)=g_j$ for each $g_j\in G_j$, $\pi ^k_l\circ \pi ^j_k=\pi
^j_l$ for each $l<k<j$ in $\Omega $. This arise questions for a
subsequent research. Whether a nonassociative analog of a
topological group has this property or not? How weak may a
nonassociative structure be that do not satisfy this property? This
article answers these questions. In it analogs of topological groups
are scrutinized with a rather mild nonassociative metagroup
structure.
\par In this article topologies on metagroups
are studied. They have specific features in comparison with groups
because of nonassociativity in general of metagroups. Homomorphisms
and quotient maps on them are investigated. In particular,
topologies on smashed products and smashed twisted wreath products
of metagroups are scrutinized in this article, which are making them
topological metagroups.  Structure of topological metagroups is
scrutinized. Moreover, their inverse homomorphism systems are
studied.
\par All main results of this paper are obtained for the first time.
They can be used for further studies of topological metagroups,
noncommutative mathematical analysis and function theory, structure
of the generalized operator algebras and $C^*$ algebras over
octonions, Cayley-Dickson algebras, operator theory and spectral
theory over octonions and Cayley-Dickson algebras, PDE,
noncommutative geometry, mathematical physics, their applications in
the sciences.

\section{Topologies on metagroups and transversal sets.}
\par {\bf Definition 2.1.} Let $G$ be a quasigroup and let $H$ be its
subquasigroup. Let $V=V_{G,H}$ be a subset in $G$ such that
\par $(2.1.1)$ $G= \bigcup_{v\in V} Hv$ and
\par $(2.1.2)$ $(Hv_1)\cap (Hv_2)=\emptyset $ for each
$v_1\ne v_2$ in $V$.
\par Then $V$ is called a transversal set of $H$ in $G$. It also
will be denoted by $D/_cA$, $V=D/_cA$ (using such notation in order
to distinguish it with $Div_r(a/b)=a/b$).

\par {\bf Remark 2.2.} For a metagroup $D$ and a submetagroup $A$
there exists the transversal set $V=V_{D,A}$ by Corollary 1 in
\cite{ludstwpmgax19}. According to Formulas $(53)$ in
\cite{ludstwpmgax19} there exist single-valued surjective maps
\par $(2.2.1)$ $\psi ^D_A: D\to A$ and
\par $(2.2.2)$ $\tau ^D_A: D\to V_{D,A}$
\\ such that $\psi ^D_A(d)=a$, $\tau ^D_A(d)=v$, where $d\in D$, $d=av$,
$a\in A$, $v\in V$. It also is denoted by $d^{\psi ^D_A}=\psi
^D_A(d)$ and $d^{\tau ^D_A}=\tau ^D_A(d)$ or shortly $d^{\psi }=\psi
(d)$ and $d^{\tau }=\tau (d)$ respectively, if $D$ and $A$ are
specified.
\par Let ${\cal C}_A(D)={\cal C}(D)\cap A$ and let ${\cal C}_1$ be
a subgroup in ${\cal C}(D)$ such that ${\cal C}_m(D)\subset {\cal
C}_1$, where ${\cal C}_m(D)$ denotes a minimal subgroup in $D$ such
that $t_D(a,b,c)\in {\cal C}_m(D)$ for each $a$, $b$ and $c$ in $D$.
For a topological metagroup $D$ using the (joint) continuity of
multiplication, $Div_l$ and $Div_r$ on $D$ one can consider without
loss of generality that ${\cal C}_1$ is closed in $D$.  This and
Condition $(1.1.9)$ imply that $A{\cal C}_1$ is a submetagroup in
$D$, since ${\cal C}_1$ is the commutative (Abelian) group and
${\cal C}_m(D)\subset {\cal C}_1$, where $PB=\{ c\in D: \exists p\in
P, \exists b\in B, c=pb \} $ for subsets $B$ and $P$ in $D$. Hence
$D/_c{\cal C}_1$ and $(A{\cal C}_1)/_c{\cal C}_1$ are the quotient
groups by Theorem 1 in \cite{ludstwpmgax19} (see also Definition 2.1
above). Therefore
\par $(2.2.3)$ $\psi ^D_A(d)=\psi ^{A{\cal C}_1}_A\circ
\psi ^D_{A{\cal C}_1}(d)$ and
\par $(2.2.4)$ $\tau ^D_A(d)=\tau ^{A{\cal C}_1}_A(\psi ^D_{A{\cal C}_1}
(d))\tau ^D_{A{\cal C}_1}(d)$
\\ for each $d\in D$. Moreover, the transversal set $V_{A{\cal C}_1,A}$ can
be chosen such that $V_{A{\cal C}_1,A}\subset V_{D,A}$ (see Remark 3
in \cite{ludstwpmgax19}).

\par {\bf Theorem 2.3.} {\it Let $G$ be a topological quasigroup with
a topology ${\cal T}_G$, let $G_1$, $G_2$ be subquasigroups such
that $G_2\subset G_1\subset G$ and let $V_1$, $V_2$, $V_{1,2}$ be
transversal sets of $G_1$, $G_2$ in $G$ and of $G_2$ in $G_1$
respectively. Then there exists a coarsest topology ${\cal
T}_{G_,G_1,G_2}$ on $G$ such that ${\cal T}_{G,G_1,G_2}\supset {\cal
T}_G$; $(G, {\cal T}_{G,G_1,G_2})$ is the topological quasigroup and
the maps $\psi ^G_{G_1}$, $\psi ^G_{G_2}$, $\psi ^{G_1}_{G_2}$,
$\tau ^G_{G_1}$, $\tau ^G_{G_2}$, $\tau ^{G_1}_{G_2}$ are
continuous, where $G$ is supplied with the topology ${\cal
T}_{G,G_1,G_2}$; $G_1$, $G_2$, $V_1$, $V_2$, $V_{1,2}$ are
considered in topologies inherited from $(G,{\cal T}_{G,G_1,G_2})$.}
\par {\bf Proof.} Note that the maps $\psi ^G_{G_1}: G\to G_1$,
$\tau ^G_{G_1}: G\to V_1$ and similarly $\psi ^G_{G_2}$, $\tau
^G_{G_2}$, $\psi ^{G_1}_{G_2}$, $\tau ^{G_1}_{G_2}$ are surjective.
\par If $X$ is a topological space with a topology ${\cal T}_X$, $f:
X\to Y$ is a surjective (single-valued) map, then it induces an
equivalence relation $\Xi _f$ on $X$ such that $x\Xi _fz$ if and
only if $f(x)=f(z)$, where $x$ and $z$ belong to $X$. Therefore,
there exists a bijection $g_Y: X/\Xi _f\to Y$, where $X/\Xi _f$
denotes the quotient topological space obtained from $X$ with the
equivalence relation $\Xi _f$. This means that $g_Y$ and the
quotient topology ${\cal T}_{X,f}$ on $X/\Xi _f$ induce the
corresponding (quotient) topology ${\cal T}^f_Y({\cal T}_X)={\cal
T}^f_Y$ on $Y$ by Proposition 2.4.3 \cite{eng}. Therefore, $f: (X,
{\cal T}_X)\to (Y, {\cal T}^f_Y)$ is continuous. More generally, if
$f_j: X\to Y$ is a surjective map for each $j\in \Lambda $, where
$\Lambda $ is a set, then there exists a coarsest topology ${\cal
T}^{(f_j: j\in \Lambda )}_Y({\cal T}_X)= {\cal T}^{(f_j: j\in
\Lambda )}_Y$ such that $f_k: (X, {\cal T}_X)\to (Y, {\cal T}^{(f_j:
j\in \Lambda )}_Y)$ is continuous and ${\cal T}^{f_k}_Y\subset {\cal
T}^{(f_j: j\in \Lambda )}_Y$ for each $k\in \Lambda $ by Proposition
2.4.14 \cite{eng}.
\par Using this we construct a new topology on $G$ by induction.
At first we take the topologies $Q^G_{G_1,1}={\cal T}^{\psi
^G_{G_1}}_{G_1}({\cal T}_{G,1})$; $Q^G_{V_1,1}={\cal T}^{\tau
^G_{G_1}}_{V_1}({\cal T}_{G,1})$; $Q^G_{G_2,1}={\cal T}^{\psi
^G_{G_2}}_{G_2}({\cal T}_{G,1})$; $Q^G_{V_2,1}={\cal T}^{\tau
^G_{G_2}}_{V_2}({\cal T}_{G,1})$; $Q^{G_1}_{G_2,1}={\cal T}^{\psi
^{G_1}_{G_2}}_{G_2}({\cal T}_{G_1,0})$; $Q^{G_1}_{V_{1,2},1}={\cal
T}^{\tau ^{G_1}_{G_2}}_{V_{1,2}}({\cal T}_{G_1,0})$; where ${\cal
T}_{G_1,0}={\cal T}_{G,1}\cap G_1$; ${\cal T}_{G_2,0}={\cal
T}_{G,1}\cap G_2$; ${\cal T}_{G,1}={\cal T}_G$.
 We put ${\cal T'}_{G_1,1}$, ${\cal T}_{V_1,1}$, ${\cal T}_{G_2,1}$,
 ${\cal T}_{V_2,1}$, ${\cal T}_{V_{1,2},1}$ to be the coarsest
 topologies generated by $Q^G_{G_1,1}\cup {\cal T}_{G_1,0}$;
 $Q^G_{V_1,1}\cup ({\cal T}_{G,1}\cap V_1)$; $Q^G_{G_2,1}\cup
 Q^{G_1}_{G_2,1}\cup {\cal T}_{G_2,0}$; $Q^G_{V_2,1}\cup ({\cal
 T}_{G,1}\cap V_2)$; $Q^{G_1}_{V_{1,2},1}\cup ({\cal T}_{G_1,0}\cap
 V_{1,2})$ respectively. Then we consider $f_1=m_G$, $f_2=Div_l$, $f_3=Div_r$,
 where $m_G$ denotes multiplication on $G$, $m_G(a,b)=ab$ for
 each $a$ and $b$ in $G$, also $f_1: G_1\times V_1\to G$,
 $f_1: G_2\times V_2\to G$, $f_1: G_2\times V_{1,2}\to G_1$, etc,
 where $G_1\times V_1$, $G_2\times V_2$, $G_2\times V_{1,2}$ are in
 the Tychonoff product topologies ${\cal T}_{G_1\times V_1,1}$;
 ${\cal T}_{G_2\times V_2,1}$; ${\cal T}_{G_2\times V_{1,2},1}$ of
 the topological spaces $(G_1,{\cal T'}_{G_1,1})\times (V_1,{\cal
 T}_{V_1,1})$; $(G_2,{\cal T}_{G_2,1})\times (V_2,{\cal
 T}_{V_{2,1}})$; $(G_2,{\cal T}_{G_2,1})\times (V_{1,2}, {\cal
 T}_{V_{1,2},1})$ respectively. This induces the topologies ${\cal
 T} ^{(f_1,f_2,f_3)}_G({\cal T}_{G_1\times V_1,1})$; ${\cal
 T} ^{(f_1,f_2,f_3)}_G({\cal T}_{G_2\times V_2,1})$; ${\cal
 T} ^{(f_1,f_2,f_3)}_{G_1}({\cal T}_{G_2\times V_{1,2},1})$ on $G$ and
 $G_1$ respectively and hence the coarsest topology $Q_{G,1}$ on $G$
 such that $Q_{G,1}\supset {\cal T}_{G,1}\cup {\cal
 T}^{(f_1,f_2,f_3)}_G({\cal T}_{G_1\times V_1,1})\cup {\cal
 T}^{(f_1,f_2,f_3)}_G({\cal T}_{G_2\times V_2,1})$, consequently,
 the coarsest topology ${\cal T}_{G,2}$ on $G$ such that ${\cal
 T}_{G,2}\supset {\cal T}_{G,1}\cup {\cal
 T}^{(f_1,f_2,f_3)}_G(Q_{G\times G,1})$, where $Q_{G\times G,1}$ is
 the Tychonoff product topology on $G\times G$ induced by $Q_{G,1}$.
\par Then there exists the coarsest topology ${Q'}_{G_1,1}$ on $G_1$ such that
${Q'}_{G_1,1}\supset {\cal T'}_{G_1,1}\cup {\cal
T}^{(f_1,f_2,f_3)}_{G_1}({\cal T}_{G_2\times V_{1,2},1})$ and hence
the coarsest topology ${\cal T}_{G_1,1}$ such that ${\cal
T}_{G_1,1}\supset {\cal T'}_{G_1,1}\cup {\cal
T}^{(f_1,f_2,f_3)}_{G_1}({Q'}_{G_1\times G_1,1})$, where
${Q'}_{G_1\times G_1,1}$ is the Tychonoff product topology on
$G_1\times G_1$ induces by ${Q'}_{G_1,1}$.
\par Assume that on the $n$-th step the topologies ${\cal T}_{G,n}$,
${\cal T}_{G_1,n-1}$, ${\cal T}_{G_2,n-1}$ on $G$, $G_1$, $G_2$
respectively are constructed, where $n\ge 2$ is the natural number.
Then we take the topologies $Q^G_{G_1,n}={\cal T}^{\psi
^G_{G_1}}_{G_1}({\cal T}_{G,n})$; $Q^G_{V_1,n}={\cal T}^{\tau
^G_{G_1}}_{V_1}({\cal T}_{G,n})$; $Q^G_{G_2,n}={\cal T}^{\psi
^G_{G_2}}_{G_2}({\cal T}_{G,n})$; $Q^G_{V_2,n}={\cal T}^{\tau
^G_{G_2}}_{V_2}({\cal T}_{G,n})$; $Q^{G_1}_{G_2,n}={\cal T}^{\psi
^{G_1}_{G_2}}_{G_2}({\cal T '}_{G_1,n-1})$;
$Q^{G_1}_{V_{1,2},n}={\cal T}^{\tau ^{G_1}_{G_2}}_{V_{1,2}}({\cal T
'}_{G_1,n-1})$; where ${\cal T '}_{G_1,n-1}$ is the coarsest
topology on $G_1$ generated by ${\cal T}_{G_1,n-1}\cup ({\cal
T}_{G,n}\cap G_1)$. Then we put ${\cal T''}_{G_1,n}$, ${\cal
T}_{V_1,n}$, ${\cal T}_{G_2,n}$,
 ${\cal T}_{V_2,n}$, ${\cal T}_{V_{1,2},n}$ to be the coarsest
 topologies generated by $Q^G_{G_1,n}\cup {\cal T '}_{G_1,n-1}$;
 $Q^G_{V_1,n}\cup ({\cal T}_{G,n}\cap V_1)$; $Q^G_{G_2,n}\cup
 Q^{G_1}_{G_2,n}\cup ({\cal T '}_{G_1,n-1}\cap G_2)\cup {\cal T}_{G_2,n-1}$;
 $Q^G_{V_2,n}\cup ({\cal  T}_{G,n}\cap V_2)$; $Q^{G_1}_{V_{1,2},n}\cup ({\cal T '}_{G_1,n-1}\cap
 V_{1,2})$ respectively. Using these topologies we get the topologies
${\cal T} ^{(f_1,f_2,f_3)}_G({\cal T}_{G_1\times V_1,n})$; ${\cal
 T} ^{(f_1,f_2,f_3)}_G({\cal T}_{G_2\times V_2,n})$; ${\cal
 T} ^{(f_1,f_2,f_3)}_{G_1}({\cal T}_{G_2\times V_{1,2},n})$ on $G$ and
 $G_1$ respectively, consequently, the coarsest topology $Q_{G,n}$ on $G$
 such that $Q_{G,n}\supset {\cal T}_{G,n}\cup {\cal
 T}^{(f_1,f_2,f_3)}_G({\cal T}_{G_1\times V_1,n})\cup {\cal
 T}^{(f_1,f_2,f_3)}_G({\cal T}_{G_2\times V_2,n})$,
where ${\cal T}_{G_1\times V_1,n}$ is the Tychonoff product topology
on the topological space $(G_1, {\cal T ''}_{G_1,n})\times (V_1,
{\cal T}_{V_1,n})$; ${\cal T}_{G_2\times V_2,n}$ corresponds to
$(G_2, {\cal T}_{G_2,n})\times (V_2, {\cal T}_{V_2,n})$, etc. This
induces the coarsest topology ${\cal T}_{G,n+1}$ on $G$ such that
${\cal  T}_{G,n+1}\supset {\cal T}_{G,n}\cup {\cal
 T}^{(f_1,f_2,f_3)}_G(Q_{G\times G,n})$, where $Q_{G\times G,n}$ is
 the Tychonoff product topology on $G\times G$ induced by $Q_{G,n}$.
\par Similarly let ${Q'}_{G_1,n}$ be the coarsest topology on $G_1$ such that
${Q'}_{G_1,n}\supset {\cal T''}_{G_1,n}\cup {\cal
T}^{(f_1,f_2,f_3)}_{G_1}({\cal T}_{G_2\times V_{1,2},n})$ and let
${\cal T}_{G_1,n}$ be the coarsest topology such that ${\cal
T}_{G_1,n}\supset {\cal T''}_{G_1,n}\cup {\cal
T}^{(f_1,f_2,f_3)}_{G_1}({Q'}_{G_1\times G_1,n})$. This process is
repeated by induction in $n\in {\bf N}$. Therefore ${\cal
T}_{G,n+1}\supset {\cal T}_{G,n}$; ${\cal T}_{G_j,n}\supset {\cal
T}_{G_j,n-1}$, ${\cal T}_{V_j,n+1}\supset {\cal T}_{V_j,n}$ and
${\cal T}_{V_{1,2},n+1}\supset {\cal T}_{V_{1,2},n}$ for each $n\in
{\bf N}$ and $j\in \{ 1, 2 \} $. Hence there exists the coarsest
topology ${\cal T}_{G,\infty }$ on $G$ such that ${\cal T}_{G,\infty
}\supset \bigcup_{n=1}^{\infty }{\cal T}_{G,n}.$ This implies that
the coarsest topology ${\cal T}_{G,G_1,G_2}$ exists such that ${\cal
T}_{G,\infty }\supset {\cal T}_{G,G_1,G_2}\supset {\cal T}_G$
relative to which the maps $f_i: G\times G\to G$ are continuous for
each $i\in \{ 1, 2 ,3 \} $ and $\psi ^G_{G_1}: G\to G_1$, $\psi
^G_{G_2}: G\to G_2$, $\tau ^G_{G_1}: G\to V_1$; $\tau ^G_{G_2}: G\to
V_2$; $\tau ^{G_1}_{G_2}: G_1\to V_{1,2}$ are continuous, where
$G_1$, $G_2$, $V_1$, $V_2$, $V_{1,2}$ are considered in the
topologies inherited from $(G, {\cal T}_{G,G_1,G_2})$.

\par {\bf Corollary 2.4.} {\it Let $D$ be a $T_1$ topological
metagroup with a topology ${\cal T}_D$ and let the conditions of
Remark 2.2 be satisfied. Then there exists the coarsest topology
${\cal T}_{D,A{\cal C}_1,A}$ on $D$ such that ${\cal T}_{D,A{\cal
C}_1,A}\supset {\cal T}_D$, $(D, {\cal T}_{D,A{\cal C}_1,A})$ is the
topological metagroup and the maps $m_D$, $Div_l$, $Div_r$, $\psi
^D_A$, $\tau ^D_A$, $\psi ^D_{A{\cal C}_1}$, $\tau ^D_{A{\cal
C}_1}$, $\psi ^{A{\cal C}_1}_A$, $\tau ^{A{\cal C}_1}_A$ are
(jointly) continuous, where $A$, $A{\cal C}_1$, $V_{D,A{\cal C}_1}$,
$V_{D,A}$, $V_{A{\cal C}_1,A}$ are considered in the topologies
inherited from $(D, {\cal T}_{D,A{\cal C}_1,A})$.}
\par {\bf Proof.} If in Theorem 2.3 $(G, {\cal T}_G)$ is a
topological metagroup, then from Formula $(1.1.9)$ it follows that
$t_G: G^3\to G$ also is (jointly) continuous, consequently, $(G,
{\cal T}_{G,G_1,G_2})$ is the topological metagroup. Therefore, the
assertion of this corollary follows from Theorem 2.3 by taking
$G=D$, $G_1=A{\cal C}_1$, $G_2=A$, $V_1=V_{D,A{\cal C}_1}$,
$V_2=V_{D,A}$, $V_{1,2}=V_{A{\cal C}_1,A}$.

\par {\bf Corollary 2.5.} {\it Let $G$ be a topological
quasigroup and let $A$ be a subquasigroup in $G$ and $V_{G,A}$ be a
transversal set of $A$ in $G$. Then there exists the coarsest
topology ${\cal T}_{G,A}$ on $G$ such that the maps $m_G$, $Div_l$,
$Div_r$, $\psi ^G_A$, $\tau ^G_A$ are continuous, where $A$,
$V_{G,A}$ are considered in the topologies inherited from $(G, {\cal
T}_{G,A})$.}
\par The assertion of this corollary is the particular case of
Theorem 2.3 with $G_2=G_1=A$ and ${\cal T}_{G,A}={\cal
T}_{G,G_1,G_2}$.

\par {\bf Remark 2.6.}  In view of Corollary 1 in
\cite{ludstwpmgax19} and Remark 2.2
\par $(2.6.1)$ $\psi ^D_A(\psi ^D_A(d))=\psi ^D_A(d)$
and $\tau ^D_A(\tau ^D_A(d))=\tau ^D_A(d)$ and
\par $(2.6.2)$ $d=d^{\psi ^D_A}d^{\tau ^D_A}$ for each $d\in D$.
\\ In particular, $(d^{\tau ^D_A})^{\psi ^D_A}=e$ and $(d^{\psi
^D_A})^{\tau ^D_A}=e$ for each $d\in D$, where $D$ is the metagroup.
Denoting $a=d^{\psi }$ and $v=d^{\tau }$, $a=a_d$, $v=v_d$, where
$\psi =\psi ^D_A$, $\tau =\tau ^D_A$, one gets $a=d/v$.
\par From Theorem 1 in \cite{ludstwpmgax19} it follows
that ${\cal C}_1^{\tau }$ is isomorphic with ${\cal C}_1/_c{\cal
C}_{1,A}$, where $B^{\tau }=\{ b^{\tau }: b\in B \} $ for $B\subset
D$. Moreover, $(A{\cal C}_1)/_c{\cal C}_{1,A}$ and $A/_c{\cal
C}_{1,A}$ are the quotient groups such that $A/_c{\cal
C}_{1,A}\hookrightarrow (A{\cal C}_1)/_c{\cal C}_{1,A}$, where
${\cal C}_{1,A}={\cal C}_1\cap A$.
\par If $a_1n_1=b(a_2n_2)$, where $a_1, a_2, b$ belong to $A$, $n_1,
n_2$ are in ${\cal C}_1$, then $(a_2\setminus (b\setminus
a_1))n_1=n_2$, hence $(a_2\setminus (b\setminus a_1))\in A\cap {\cal
C}_1={\cal C}_{1,A}$. Vive versa if $n_1=\beta n_2$ with $\beta \in
{\cal C}_{1,A}$, $n_1, n_2$ in ${\cal C}_1$, then for each $a_1,
a_2$ in $A$ there exists $b=(a_1\beta )/a_2\in A$ such that
$a_1\beta =ba_2$ and consequently, $a_1n_1=a_1\beta n_2=b(a_2n_2)$.
Thus the quotient groups $(A{\cal C}_1)/_cA$ and ${\cal C}_1/_c{\cal
C}_{1,A}$ are isomorphic. From the latter isomorphism, Remark 3 in
\cite{ludstwpmgax19} it follows that $V_{D,A}$, $V_{{\cal C}_1,{\cal
C}_{1,A}}$ and $V_{D,A{\cal C}_1}$ can be chosen such that
\par $(2.6.3)$ $V_{{\cal C}_1,{\cal C}_{1,A}}V_{D,A{\cal C}_1}=V_{D,A}$; \hspace{0.5cm}
$V_{{\cal C}_1,{\cal C}_{1,A}}=V_{A{\cal C}_1,A}$, \\ since ${\cal
C}_1$ is the invariant subgroup in $D$ and $p_D(a,b,c)=e$ if $e\in
\{ a, b , c \} \subset D$. Hence
\par $(2.6.4)$ $(d^{\psi }\gamma )^{\psi }=d^{\psi }\gamma ^{\psi }$
and $(d^{\psi }\gamma )^{\tau }=\gamma ^{\tau }$ for each $d\in D$
and $\gamma \in {\cal C}_1$, since $d^{\psi }\gamma = d^{\psi
}(\gamma ^{\psi }\gamma ^{\tau })=(d^{\psi }\gamma ^{\psi })\gamma
^{\tau }$, since $\gamma ^{\psi } \in {\cal C}_{1,A}$ and $\gamma
^{\tau }=\gamma ^{\psi }\setminus \gamma \in {\cal C}_1$. We remind
that
\par $(2.6.5)$ $(a^{\tau })^{[c]}:=(a^{\tau }c)^{\tau }=: \nu
_{D,A}(a,c)$ \\ for each $a$ and $c$ in $D$ (see Formula $(68)$ in
\cite{ludstwpmgax19}). Suppose that the conditions of Remark 4 in
\cite{ludstwpmgax19} are satisfied. Let on the Cartesian product
$C=D\times F$ (or $C^*=D\times F^*$) for each $d$, $d_1$ in $D$,
$f$, $f_1$ in $F$ (or $F^*$ respectively) a binary operation be:
\par $(2.6.6)$ $(d_1,f_1)(d,f)=(d_1d, \xi ((d_1^{\psi }, f_1),
(d^{\psi },f))f_1f^{ \{ d_1 \} })$, \\ where $\xi ((d_1^{\psi },
f_1), (d^{\psi },f))(v)= \xi ((d_1^{\psi }, f_1(v)), (d^{\psi
},f(v)))\in {\cal C}_1$ for each $v\in V_{D,A}$ (see Formula $(86)$
in \cite{ludstwpmgax19}).

\par {\bf Theorem 2.7.} {\it Let $G$ be a $T_1$ topological
quasigroup and let ${\cal B}=\bigcup_{g\in G} {\cal B}_g$, where
${\cal B}_g$ is an open base at $g$ in $G$. Then ${\cal B}$
satisfies the following properties $(2.7.1)$-$(2.7.8)$: \par
$(2.7.1)$ $\forall g\in G,$ $\forall h\in G$, ${\cal B}_{hg}=h{\cal
B}_g$ and ${\cal B}_{gh}={\cal B}_gh$; \par $(2.7.2)$ $\forall g\in
G$, $\forall h\in G$, ${\cal B}_{h\setminus g}=h\setminus {\cal
B}_g$ and ${\cal B}_{g/h}={\cal B}_g/h$; \par $(2.7.3)$ $\forall
g\in G$, $\forall B\in {\cal B}_g$, $g\in B$; \par $(2.7.4)$
$\forall g\in G$, $\forall U_g\in {\cal B}_g$, $\forall a\in G$,
$\exists U_a\in {\cal B}_a$, $\exists b=a\setminus g\in G$, $\exists
U_b\in {\cal B}_b$, $U_aU_b\subset U_g$; \par $(2.7.5)$ $\forall
g\in G$, $\forall U_g\in {\cal B}_g$, $\forall a\in G$, $\exists
U_a\in {\cal B}_a$, $\exists b=g\setminus a\in G$, $\exists U_b\in
{\cal B}_b$, $U_a/U_b\subset U_g$; \par $(2.7.6)$ $\forall g\in G$,
$\forall U_g\in {\cal B}_g$, $\forall a\in G$, $\exists U_a\in {\cal
B}_a$, $\exists b=a/ g\in G$, $\exists U_b\in {\cal B}_b$,
$U_b\setminus U_a\subset U_g$;
\par $(2.7.7)$ $\forall g\in G$, $\forall U\in {\cal B}_g$, $\forall V\in {\cal B}_g$,
$\exists W\in {\cal B}_g$, $W\subset U\cap V$; \par $(2.7.8)$
$\forall g\in G$, $\bigcap_{U\in {\cal B}_g} = \{ g \} $. \par
Conversely, let $G$ be a quasigroup and let ${\cal B}$ be a family
of subsets in $G$ satisfying $(2.7.1)$-$(2.7.8)$. Then the family
${\cal B}$ is a base for a $T_1$ topology ${\cal T}_{\cal B}={\cal
T}_{\cal B}(G)$ on $G$ and $(G,{\cal T}_{\cal B})$ is a topological
quasigroup.}
\par {\bf Proof.} Assume that $G$ is a topological quasigroup
with a topology ${\cal T}$.
For each $g\in G$ one can take
\par $(2.7.9)$ ${\cal B}_g= \{ U: U\in {\cal T}, g\in U \} $ and put
${\cal B}=\bigcup_{g\in G} {\cal B}_g$.
\par Let $L_q: G\to G$ be a left shift map, $L_qg=qg$, $R_q: G\to G$
be a right shift map, $R_qg=gq$ for each $g\in G$, where $q\in G$.
Then we put $\check{L}_qg=q\setminus g$, $\check{R}_qg=g/q$ for each
$g\in G$, where $q\in G$. Hence
\par $(2.7.10)$ $L_q(\check{L}_qg)=g$, $\check{L}_q(L_qg)=g$,
$R_q(\check{R}_qg)=g$, $\check{R}_q(R_qg)=g$ for each $g\in G$. The
maps $m_G$, $Div_l$ and $Div_r$ from $G\times G$ into $G$ are
(jointly) continuous, where $m_G$ denotes multiplication on $G$.
From $(2.7.9)$ it follows that $L_q$, $\check{L}_q$, $R_q$,
$\check{R}_q$ are homeomorphisms from $G$ onto $G$ as the
topological spaces for each $q\in G$. This implies $(2.7.1)$,
$(2.7.2)$, while $(2.7.3)$ follows from $(2.7.9)$. Then
$(2.7.4)$-$(2.7.6)$ follow from the continuity of $m_G$, $Div_l$,
$Div_r$. Certainly, $(2.7.9)$ implies $(2.7.7)$. Property $(2.7.8)$
follows from $G$ being $T_1$ as the topological space and from
${\cal B}_g$ being the open base at $g$ in $G$.
\par Vice versa, let ${\cal B}$ be a family of subsets of $G$
satisfying conditions $(2.7.1)$-$(2.7.8)$. Let ${\cal T}$ be a
family such that \par $(2.7.11)$ $\forall A\in {\cal T}$, $\exists
{\cal E}\subset {\cal B}$, $\bigcup_{U\in \cal E} U=A$.
\par Assume that ${\cal F}\subset {\cal T}$, hence $\bigcup_{U\in
\cal F}U\in {\cal T}$, since $\forall C\in {\cal F}$, $\exists {\cal
F}_C\subset {\cal F}$, $\bigcup_{U\in {\cal F}_C}U=C$. Then we take
any fixed $W_1\in {\cal T}$ and $W_2\in {\cal T}$ and put $W=W_1\cap
W_2$. For each $g\in W$ there exists $B\in {\cal B}_g$ such that
$g\in B\subset W$, consequently, there exists ${\cal E}\subset {\cal
T}$ with $\bigcup_{U\in \cal E}U=W$. Thus ${\cal T}$ is a topology
on $G$. Then for each $A$ and ${\cal E}$: $\forall U\in {\cal E}$,
$\forall g\in U$, $\exists U_g\in {\cal B}_g$, $g\in U_g\subset U$,
$\forall x\in G$, $U_gx\in {\cal B}_{gx}$, $xU_g\in {\cal B}_{xg}$,
$x\setminus U_g\in {\cal B}_{x\setminus g}$, $U_g/x\in {\cal
B}_{g/x}$ by $(2.7.1)$, $(2.7.2)$ and $(2.7.8)$, similarly to
$(2.7.11)$.
\par Thus the family ${\cal B}$ is a base for the topology ${\cal
T}$. From $(2.7.11)$ and $(2.7.4)$-$(2.7.6)$ we infer that $m_G$,
$Div_l$, $Div_r$ are (jointly) continuous maps with respect to
${\cal T}$. Then $(2.7.11)$ and $(2.7.1)$, $(2.7.2)$ imply that
$\forall g\in G$, $\forall h\in G$, ${\cal T}_{hg}=h{\cal T}_g$,
${\cal T}_{gh}={\cal T}_gh$, ${\cal T}_{h\setminus g} =h\setminus
{\cal T}_g$, ${\cal T}_{g/h}={\cal T}_g/h$, where ${\cal T}_g = \{
A\in {\cal T}: g\in A \} $. Thus $(G,{\cal T})$ is the topological
quasigroup. From $(2.7.11)$ and $(2.7.8)$ it follows that $(G, {\cal
T})$ is $T_1$ as the topological space.

\par {\bf Lemma 2.8.} {\it Let $B$ be a quasigroup, let $V$ be a set,
$F=B^V = \{ f: V\to B \} $ be a family of all (single-valued) maps
from $V$ into $B$, $W(S,Q) := \{ f\in F: f(S)\subset Q \} $, where
$V\supset S\ne \emptyset $, $B\supset Q\ne \emptyset $. Then for
nonvoid subsets $S$, $S_1$, $S_2$, $S_i$ in $V$, $Q$, $Q_1$, $Q_i$
in $B$:
\par $(2.8.1)$ $\forall S\subset V$, $\forall Q_2\subset B$, $\forall
Q_1\subset Q_2$, $W(S,Q_1)\subset W(S,Q_2)$; \par $(2.8.2)$ $\forall
b\in B$, $\forall S\subset V$, $\forall Q\subset B$,
$W(S,b/Q)=b/W(S,Q)$ $\& $ $W(S,Q\setminus b)=W(S,Q)\setminus b$;
\par $(2.8.3)$ $\forall S\subset V$, $\forall Q\subset B$, $\forall
Q _1\subset B$, $W(S,Q)\setminus W(S,Q_1)\subset W(S,Q\setminus
Q_1)$ $\& $ $W(S,Q)/W(S,Q_1)\subset W(S,Q/Q_1)$ $\& $
$W(S,Q)W(S,Q_1)\subset W(S,QQ_1)$;
\par $(2.8.4)$ $\forall S_2\subset V$, $\forall S_1\subset S_2$,
$\forall Q\subset B$, $W(S_2,Q)\subset W(S_1,Q)$;
\par $(2.8.5)$ for each set $\Lambda $, $\forall i\in \Lambda $,
$\forall S_i\subset V$, $\forall Q\subset B$, $W(\bigcup_{i\in
\Lambda } S_i, Q)=\bigcap_{i\in \Lambda } W(S_i,Q)$;
\par $(2.8.6)$ for each set $\Lambda $, $\forall i\in \Lambda $,
$\forall Q_i\subset B$, $(\bigcap_{i\in \Lambda }Q_i\ne \emptyset
)$, $\forall S\subset V$, $W(S,\bigcap_{i\in \Lambda }Q_i)
=\bigcap_{i\in \Lambda }W(S,Q_i)$.}
\par {\bf Proof.} $(2.8.1)$. From $f(S)\subset Q_1$ it follows that
$f(S)\subset Q_2$.
\par $(2.8.2)$. $(\forall f\in W(S,b/Q),$ $\forall s\in S$, $\exists
q\in Q$, $(f(s)=b/q \leftrightarrow q=f(s)\setminus b))$
$\rightarrow $ $(\exists f_1\in W(S,Q),$ $(f=b/f_1 \leftrightarrow
f_1=f\setminus b)$;
\par $(\forall f_1\in W(S,Q)$, $\forall s\in S$, $\exists q_1\in Q$,
$f_1(s)=q_1)$ $\rightarrow $ $(\exists f\in W(S,b/Q),$
$(f_1=f\setminus b \leftrightarrow f=b/f_1))$.
\par Thus $W(S,b/Q)=b/W(S,Q)$. Symmetrically it is proved
$W(S,Q\setminus b)=W(S,Q)\setminus b$.
\par $(2.8.3)$. $(\forall f\in W(S,Q),$ $\forall f_1\in W(S,Q_1),$
$\forall s\in S,$ $f(s)\setminus f_1(s)\in Q\setminus Q_1)$
$\rightarrow $ $(\forall f\in W(S,Q),$ $\forall f_1\in W(S,Q_1),$
$f\setminus f_1\in W(S,Q\setminus Q_1))$;
\par $(\forall f\in W(S,Q)$, $\forall f_1\in W(S,Q_1),$ $\forall
s\in S$, $f(s)/f_1(s)\in Q/Q_1)$ $\rightarrow $ $(\forall f\in
W(S,Q),$ $\forall f_1\in W(S,Q_1)$, $f/f_1\in W(S,Q/Q_1))$;
\par $(\forall f\in W(S,Q)$, $\forall f_1\in W(S,Q_1)$, $\forall
s\in S$, $f(s)f_1(s)\in QQ_1)$ $\rightarrow $ $(\forall f\in
W(S,Q),$ $\forall f_1\in W(S,Q_1),$ $ff_1\in W(S,QQ_1))$,
\\ since $Q\setminus Q_1= \bigcup ( b\in B: \exists q\in Q, \exists
q_1\in Q_1, b=q\setminus q_1 ) $, \par $Q/Q_1= \bigcup ( b\in B:
\exists q\in Q, \exists q_1\in Q_1, b=q/q_1 ) $,
\par $QQ_1= \bigcup ( b\in B: \exists q\in Q, \exists
q_1\in Q_1, b=qq_1 ) $.
\par $(2.8.4)$. $(\forall f\in W(S_2,Q)$, $f(S_2)\subset Q)$
$\rightarrow $ $(f(S_1)\subset Q)$,
\\ hence $W(S_2,Q)\subset W(S_1,Q)$, since $S_1\subset S_2\subset
V$.
\par $(2.8.5)$. $\forall f\in W(\bigcup_{i\in \Lambda }S_i,Q)$,
$((f(\bigcup_{i\in \Lambda }S_i)\subset Q)$ $\leftrightarrow $
$(\forall i\in \Lambda ,$ $f(S_i)\subset Q))$;
\par $\forall f\in \bigcap_{i\in \Lambda }W(S_i,Q),$ $((\forall i\in
\Lambda ,$ $f(S_i)\subset Q)$ $\leftrightarrow $ $(f(\bigcup_{i\in
\Lambda }S_i)\subset Q))$.
\par $(2.8.6)$. $\forall s\in S$, $((\forall i\in \Lambda ,$ $f(s)\in
Q_i)$ $\leftrightarrow $ $(f(s)\in \bigcap_{i\in \Lambda }Q_i))$,
\\ hence $W(S,\bigcap_{i\in \Lambda }Q_i)=\bigcap_{i\in \Lambda
}W(S,Q_i)$.

\par {\bf Remark 2.9.} Let the conditions of Theorem 5 in
\cite{ludstwpmgax19} be satisfied. We consider the topology ${\cal
T}_{D,A{\cal C}_1,A}$ on the topological metagroup $D$, also $A$,
${\cal C}_1$ and the transversal set $V=V_{D,A}$ in the topology
inherited from $(D,{\cal T}_{D,A{\cal C}_1,A})$. Let $C(V,B)$ denote
a family of all continuous maps $f: V\to B$, where $B$ is a
topological metagroup with a topology ${\cal T}_B$. Assume that
$(D,{\cal T}_{D,A{\cal C}_1,A})$ and $(B,{\cal T}_B)$ are $T_1$ as
topological spaces. \par As usually $U$ is a canonical closed subset
(i.e. a closed domain) in $V$ if and only if $U=cl_V(Int_VU)$, where
$Int_VU$ denotes the interior of $U$ in $V$, while $cl_VS$ denotes
the closure of $S$ in $V$, where $S\subset V$. Let ${\cal F}$ be a
family of nonvoid canonical closed subsets $U$ in $(V, {\cal
T}_{D,A{\cal C}_1,A}\cap V)$ such that
\par $(2.9.1)$ $\forall U_1\in {\cal F}$, $\forall U_2\in {\cal F}$,
$U_1\cup U_2\in {\cal F}$;
\par $(2.9.2)$ $\forall v\in V$, $\forall v\in P\in ({\cal T}_{D,A{\cal C}_1,A}\cap
V)$, $\exists U_1\in {\cal F}$, $v\in Int (U_1)$ $\& $ $U_1\subset
P$.
\par We put $C_t= \bigcup ((d,f)\in C: d\in D, f\in C(V,B))$, where
$C=D\Delta _A^{\phi , \eta , \kappa , \xi }F=D\Delta ^{\phi , \eta ,
\kappa , \xi }F$ (see Definition 5 in \cite{ludstwpmgax19}). Let
${\cal W}$ be a family of all subsets $W(S,Q)$ in $C(V,B)$ such that
$S\in {\cal F}$ and $Q\ne \emptyset $ is open in $(B,{\cal T}_B)$.

\par {\bf Theorem 2.10.} {\it Let the conditions of Remark 2.9 be
satisfied and let the maps $\phi , \eta , \kappa , \xi $ be jointly
continuous (see Remark 1 in \cite{ludstwpmgax19}). Then ${\cal
T}_{D,A{\cal C}_1,A}\times {\cal W}$ is a base of a topology ${\cal
T}_{C_t}$ on $C_t$ relative to which $C_t$ is a topological $T_1\cap
T_3$ loop.}
\par {\bf Proof.} Evidently each constant map $f_b: V\to B$
belongs to $C(V,B)$, where $b$ in $B$ is arbitrary fixed, $f_b(v)=b$
for each $v\in V$. This induces the natural embedding of $B$ into
$C(V,B)$, consequently, $C(V,B)$ is the nonvoid metagroup with
pointwise multiplication and $Div_l$ and $Div_r$. If $U_1$ and $U_2$
belong to ${\cal F}$, then $U_1\cup U_2\in {\cal F}$ (see \S I.8.8
in \cite{kuratb}, 1.1.C \cite{eng}). \par By the conditions of
Remark 2.9 $ \{ Int_V(U): U\in {\cal F} \} $ is a base of the
topology on $(V, {\cal T}_{D,A{\cal C}_1,A}\cap V)$. From $D$ and
$B$ being the $T_1$ topological metagroups and hence topological
quasigroups, it follows that they are $T_3$ as the topological
spaces. Therefore for each open subset $S_0$ in $D$ (or $V$) and
each $d_0\in S_0$ there exists a canonical closed subset (i.e.
closed domain) $S$ in $D$ (or $V$ respectively) such that $d_0\in
Int_DS$ (or $d_0\in Int_VS$ respectively) and $S\subset S_0$ by
Proposition 1.5.5 \cite{eng}. Hence ${\cal W}$ is a base of a
topology on $C(V,B)$. By virtue of Theorem 5 in \cite{ludstwpmgax19}
and Corollary 2.4 the maps $\tau ^D_A$ and $\psi ^D_A$ are
continuous on $D$ and $(D, {\cal T}_{D,A{\cal C}_1,A})$ is the
topological metagroup.
\par In view of Lemma 6 in \cite{ludstwpmgax19} the map
$w_j: D\times D\times V_{D,A}\to {\cal C}_1$ is (jointly) continuous
as the composition of jointly continuous maps for each $j\in \{ 1,
2, 3 \} $. According to Remark 4 in \cite{ludstwpmgax19} $f^{ \{ d
\} } (v)=f^{s(d,v)} (v^{[d\setminus e]})=\phi
(s(d,v))f(v^{[d\setminus e]})$ for each $d\in D$, $v\in V$, $f\in
F$. The map $A\times B\mapsto \phi (a)b\in B$ is jointly continuous
by the conditions of this theorem. Lemma 2.8 imply that ${\cal W}$
is the base of a topology on $C(V,B)$. We take the topology ${\cal
T}_{C_t}$ generated by the base ${\cal B}_{D,A{\cal C}_1,A}\times
{\cal W}$, where ${\cal B}_{D,A{\cal C}_1,A}$ denotes the base of
the topology ${\cal T}_{D,A{\cal C}_1,A}$.
\par Hence Lemma 2.8 and Formula $(2.6.6)$ imply that
multiplication $m_{C_t}: C_t\times C_t\to C_t$ is (jointly)
continuous. From Formulas in the proof of Theorem 5 in
\cite{ludstwpmgax19} it follows that $Div_l$ and $Div_r$ are jointly
continuous from $C_t\times C_t$ into $C_t$.
\par The base ${\cal W}$ generates the $T_1$ topology ${\cal T}_{\cal W}$
on $C(V,B)$ by Lemma 2.8 and Remark 2.9. Hence ${\cal T}_{C_t}$ is
the $T_1$ topology on $C_t$, since $(D,{\cal T}_{D,A{\cal C}_1,A})$
is $T_1$ by the conditions of this theorem. This implies that $(C_t,
{\cal T}_{C_t})$ is the $T_3$ topological loop.

\par {\bf Theorem 2.11.} {\it Let the conditions of Theorem 2.10 be
satisfied, let $(D, {\cal T}_{D,A{\cal C}_1,A})$ be locally compact
(or compact) and ${\cal F}$ be a family of all canonical closed
compact subsets in $(V, {\cal T}_{D,A{\cal C}_1,A}\cap V)$, let
$F_0$ be closed in $(C(V,B), {\cal T}_{\cal W})$, let also
$cl_{C(V,B)} F_0(v)$ be compact for each $v\in V$, where $F_0(v)= \{
f(v): f\in F_0 \} $, let for each compact subset $Z$ in $V$ the
restriction $F_0|_Z$ be evenly continuous, let $D\Delta _A^{\phi ,
\eta , \kappa , \xi }F_0=:C_0$ be a subloop in $C_t$. Then $C_0$ is
a locally compact (or compact respectively) loop.}
\par {\bf Proof.} Since $V=\psi ^{-1}(e)$ and $D$ is $T_1\cap T_3$ and locally compact,
then $V$ is closed in $D$ and hence locally compact by Theorem 3.3.8
\cite{eng}, consequently, $V$ is a $k$-space. From Lemma 2.8, Remark
2.9 and the conditions on ${\cal F}$ it follows that ${\cal W}$
induces a compact-open topology ${\cal T}_{\cal W}$ on $C(V,B)$. For
each compact subset $Z$ in $V$ and open $U_D$ in $D$ the conditions
of this theorem imply that $\bigcup \{ F_0^{ \{ d \} }: d\in U_D \}
|_Z$ is evenly continuous, since $F_0^{ \{ d \} }\subset F_0$ for
each $d\in D$. By virtue of Theorem 3.4.21 \cite{eng} $F_0$ is
compact. Since $(D, {\cal T}_{D,A{\cal C}_1,A})$ is locally compact,
it is sufficient to take any open $U_D$ in $D$ with the compact
closure $cl_D(U_D)$ in the ${\cal T}_{D,A{\cal C}_1,A}$ topology.
Moreover, $cl_{C(V,B)}(\bigcup \{ F_0^{ \{ d \} }(v): d\in U_D \} )$
is closed in $F_0$ and hence compact for each $v\in V$. From $B$
being $T_1$ it follows that $B$ is $T_3$.
\par From the compactness of $F_0$ and Theorem 2.10 it follows that
$C_0$ is either the locally compact loop, if $D$ is locally compact,
or the compact loop if $D$ is compact.

\par {\bf Proposition 2.12.} {\it Let the conditions of Theorem 5 in
\cite{ludstwpmgax19} be satisfied and let $i: D\to D$ and $j: B\to
B$ be automorphisms of the metagroups $D$ and $B$ such that
$i|_{{\cal C}_1}=j|_{{\cal C}_1}$. Then there exists a loop
$C_{i,j}$ and an isomorphism
\par $(2.12.1)$ $\theta _{i,j}: C\to C_{i,j}$ of $C$ onto $C_{i,j}$
such that $\theta _{i,j}|_D=i$ and $\theta _{i,j}|_B=j$.}
\par {\bf Proof.} By the conditions of this proposition
$i(ab)=i(a)i(b)$, $i(a/b)=i(a)/i(b)$ and $i(a\setminus
b)=i(a)\setminus i(b)$ for each $a$ and $b$ in $D$, similarly for
$i^{-1}$. Therefore, $i: A\to i(A)$ is an isomorphism of metagroups
and $i: {\cal C}_1\to i({\cal C}_1)$ is an isomorphism of groups
such that $i({\cal C}_1)=j({\cal C}_1)$. In view of Corollary 1 in
\cite{ludstwpmgax19} $V_{D,i(A)}=V_{i(D),i(A)}=i(V_{D,A})$,
$V_{i({\cal C}_1),i({\cal C}_{1,A})}=i(V_{{\cal C}_1,{\cal
C}_{1,A}})$, $V_{D,i(A{\cal C}_1)}=i(V_{D,A{\cal C}_1})$, $i({\cal
C}_{1,A})=i({\cal C}_1)\cap i(A)$.
\par We put
\par $(2.12.2)$ $C_{i,j}=D\Delta _A^{\phi _{i,j}, \eta _{i,j},
\kappa _{i,j}, \xi _{i,j}}F_{i,j}$ with \par $F_{i,j}= \{
j(f(i(v))): v\in V_{D,A}, f\in F \} $, $\phi _{i,j}=j\circ \phi
\circ i^{-1}$, $\phi _{i,j}: i(A)\to {\cal A}(B)$, \par $\eta
_{i,j}: i(A)\times i(A)\times B\to i({\cal C})$, $ ~ \kappa _{i,j}:
i(A)\times B\times B\to i({\cal C})$, \par $\xi _{i,j}: (i(A)\times
B)\times (i(A)\times B)\to i({\cal C})$ such that \par $\eta
_{i,j}(v',u',b')=i(\eta (i^{-1}(v'),i^{-1}(u'),j^{-1}(b')))$,
\par $\kappa _{i,j}(u',c',b')=i(\kappa
(i^{-1}(u'),j^{-1}(c'),j^{-1}(b')))$, \par $\xi
_{i,j}((u',c'),(v',b'))=i(\xi ((i^{-1}(u'),j^{-1}(c')),
(i^{-1}(v'),j^{-1}(b'))))$ for each $u'$, $v'$ in $i(A)$, $b'$, $c'$
in $B$; $\psi _i: D\to i(A)$, $\tau _i: D\to V_{D,i(A)}$ such that
\par $\psi _i(d')=i(\psi (i^{-1}(d')))$, $ ~ \tau _i(d')=i(\tau
(i^{-1}(d')))$ for each $d'\in D$. Let
\par $(2.12.3)$ $\theta _{i,j}((d,f)(v))=(i(d),j(f))(i(v))$ for each
$(d,f)\in C$ and $v\in V_{D,A}$, where $d\in D$, $f\in F$.
Therefore, $(2.6.6)$ and $(2.12.3)$ imply that
\par $(\theta _{i,j} (d_1,f_1(v))) (\theta _{i,j} (d,f(v)))
=$\par $(i(d_1d), \xi _{i,j}(((i(d_1))^{\psi _i}, j(f_1)),
((i(d))^{\psi _i}, j(f))) j(f_1) (j(f))^{ \{ i(d_1) \} } )(i(v))$
and
\par $\theta _{i,j}((d_1,f_1(v)))(d,f(v)))= $\par $(i(d_1d), i(\xi
((d_1^{\psi }, f_1), (d^{\psi },f))) j(f_1) j(f^{ \{ i(d_1) \} }
))(i(v))$ \\ for each $v\in V_{D,A}$, consequently,
\par $(2.12.4)$ $\theta _{i,j}((d_1,f_1) (d,f)) = (\theta
_{i,j}(d_1,f_1)) (\theta _{i,j}(d,f))$ \\ for each $(d,f)=g$ and
$(d_1,f_1)=g_1$ in $C$. Since $C$ and $C_{i,j}$ are loops, then
$(2.12.2)$ and $(2.12.4)$ imply that $\theta _{i,j}(g/g_1)=\theta
_{i,j}(g)/\theta _{i,j}(g_1)$ and $\theta _{i,j}(g\setminus
g_1)=\theta _{i,j}(g)\setminus \theta _{i,j}(g_1)$ for each $g$ and
$g_1$ in $C$. Thus $\theta _{i,j} : C\to C_{i,j}$ is the isomorphism
of these loops.

\par {\bf Corollary 2.13.} {\it Assume that the conditions of
Proposition 2.12 are satisfied and $i$, $i^{-1}$, $j$, $j^{-1}$ are
continuous relative to ${\cal T}_{D,A{\cal C}_1,A}$ and ${\cal T}_B$
topologies on $D$ and $B$ respectively. Then $\theta _{i,j}$ and
$\theta _{i,j}^{-1}$ are continuous relative to ${\cal T}_{C_t}$ and
${\cal T}_{C_{i,j,t}}$ topologies on $C_t$ and $C_{i,j,t}$
respectively.}
\par {\bf Proof.} This follows from Proposition 2.12, Formula $(2.12.3)$
and Theorem 2.10.
\par {\bf Remark 2.14.} If the conditions of Theorem 6
instead of that of Theorem 5 in \cite{ludstwpmgax19} are satisfied,
then in Theorems 2.10 and 2.11, Corollary 2.13 $C_t$, $C_0$,
$C_{i,j,t}$ are topological metagroups; in Proposition 2.12 $C$ and
$C_{i,j}$ are metagroups.

\section{Homomorphisms of topological metagroups.}
\par {\bf Definition 3.1.} A topological space $X$ is called homogeneous, if
for each $x\in X$ and each $y\in X$, there exists a homeomorphism
$f$ of $X$ onto itself such that $f(x)=y$.
\par If $G$ is a metagroup and $H$ is a submetagroup in $G$, $b\in
G$, then $bH$ (or $Hb$) is called a left coset (or a right coset
respectively).

\par {\bf Lemma 3.2.} {\it Let $G$ be a quasigroup and $H$ be a
subquasigroup in $G$ such that
\par $(3.2.1)$ $(Hb)a=H(ba)$ for each $a$ and $b$ in $G$.
\\ Then there exists a set $G/_cH$ of all right cosets and a
single-valued quotient map $\pi : G\to G/_cH$.}
\par {\bf Proof.} For each $a$, $b$ in $G$: \par $((Ha)\cap (Hb)\ne
\emptyset ) \leftrightarrow (\exists h_1\in H, \exists h_2\in H,
h_1a=h_2b) \leftrightarrow (\exists h_1\in H, \exists h_2\in H,
(h_1a)/b=h_2).$
\par Notice that Condition $(3.2.1)$ implies $Hb=(H(ba))/a$ for each
$a$ and $b$ in $G$, hence \par $(3.2.2)$ $H(c/a)=(Hc)/a$ for each
$a$ and $c$ in $G$, \\ since for each $a$ and $c$ in $G$ there
exists $b\in G$ with $ba=c$, that is, $b=c/a$. Therefore, $\forall
a\in G, \forall b\in G, \forall h_1\in H, \exists h_3\in H,
(h_1a)/b=h_3(a/b)$, hence
\par $\forall a\in G, \forall b\in G, (((Ha)\cap (Hb)\ne \emptyset )
\leftrightarrow (H(b/a)=H)\leftrightarrow (Ha=Hb)),$ \\ since $HH=H$
and $H(Hb)=Hb$ for each $b\in G$ by $(3.2.1)$ and $(1.1.1)$. This
implies that a set $G/_cH= \{ Hb: b\in G \} $ of all right cosets
$Hb$ exists and there exists a single-valued map $\pi : G\to G/_cH$
such that $\pi (b)=Hb$ for each $b\in G$.

\par {\bf Example 3.3.} If the conditions of Corollary 1 in
\cite{ludstwpmgax19} are satisfied, either $(G=D, H=A{\cal C}_1)$ or
$(G=A{\cal C}_1, H=A)$, then Condition $(3.2.1)$ is satisfied for
these pairs.

\par {\bf Theorem 3.4.} {\it Assume that $G$ is a topological
$T_1$ quasigroup with a topology ${\cal T}_G$. Assume also that $H$
is a closed subquasigroup satisfying Condition $(3.2.1)$ in $G$ and
$G/_cH$ is supplied with the quotient topology with respect to the
quotient map $\pi : G\to G/_cH$ (see Lemma 3.2). Then for each
arbitrary fixed $b\in G$ and $g\in G$ the family $ \{ \pi (Vg): V\in
{\cal T}_G, b\in V \} $ is a local base of $G/_cH$ at $H(bg)\in
G/_cH$. Moreover, the map $\pi $ is open and $G/_cH$ is a
homogeneous $T_1$-space.}
\par {\bf Proof.} In view of Lemma 3.2 $\pi ^{-1}(\pi (H(Vg)))=
H(Vg)=(HV)g$ for each $V\in {\cal T}_G$, consequently, $\pi
(H(Vg))=\pi ((HV)g)$ is open in $G/_cH$, since $\pi $ is the
quotient map and $(G, {\cal T}_G)$ is the topological quasigroup.
Consider any fixed $g$, $b$ in $G$ and an open neighborhood $P$ of
$H(bg)$ in $G/_cH$. Then $\pi ^{-1}(P)=:Q$ is open and $H(bg)\subset
Q$, since $\pi $ is continuous and $H(Hg)=Hg$. For each $V\in {\cal
T}_G$ such that $b\in V$ and $Vg\subset Q$ one gets $\pi (Vg)\subset
P$, consequently, $\pi ^{-1}(\pi (Vg))\subset Q$ and hence $\pi
(H(Vg))\subset P$, since $H(Vg)=\pi ^{-1}(\pi (Vg))$. Thus the map
$\pi $ is open and $ \{ \pi (Vg): V\in {\cal T}_G, b\in V \} $ is
the local base of $G/_cH$ at $H(bg)$ in $G/_cH$.  Then we put
$S_b(Hg)=(Hg)b$ for each $b\in G$, where $g\in G$. Hence $S_b:
G/_cH\to G/_cH$, since $(Hg)b=H(gb)$ by $(3.2.1)$. Therefore,
$\check{S}_bS_b=I$ and $S_b\check{S}_b=I$ for each $b\in G$, where
$\check{S}_b(Hg)=(Hg)/b$, $Ig=g$ for each $g\in G$. \par This
implies that $S_b$ is a bijection for each $b\in G$. For each $a$,
$b$ in $G$ and an open neighborhood $V$ of $g$, where $g\in G$, $\pi
(H((Va)b))$ is a basic neighborhood of $H((ga)b)=(H(ga))b$ and $\pi
(H(Va))$ is a basic neighborhood of $H(ga)$ in $G/_cH$. Therefore,
$S_b: G/_cH\to G/_cH$ is a homeomorphism, since $S_b\pi (H(Va))=\pi
((H(Va))b)=\pi (H((Va)b))$.
\par For each $a$, $b$ in $G$ we deduce from $(3.2.1)$ that
$S_b(H(a/b))=Ha$, consequently, the quotient space $G/_cH$ is
homogeneous. The space $G/_cH$ is $T_1$, since $Hb$ is closed in $G$
and $\pi $ is the quotient map.

\par {\bf Definition 3.5.} The space $G/_cH$ in Theorem 3.4 is called
the right coset space, $S_b$ is called the right translation of
$G/_cH$ by $b$.

\par {\bf Corollary 3.6.} {\it Let $G$ be a topological quasigroup
and let $H$ be a closed subquasigroup in $G$ satisfying Condition
$(3.2.1)$. Then $R_b$ and $S_b$ are homeomorphisms of $G$ and
$G/_cH$ respectively and \par $(3.6.1)$ $\pi \circ R_b=S_b\circ \pi
$ for each $b\in G$.}

\par {\bf Corollary 3.7.} {\it Assume that $G$ is a topological
quasigroup and $H$ is a closed invariant subquasigroup in $G$. Then
$G/_cH$ with the quotient topology, multiplication, $Div_l$ and
$Div_r$ induced from $G$, is a topological quasigroup and the
quotient map $\pi : G\to G/_cH$ is an open continuous homomorphism.}

\par {\bf Corollary 3.8.} {\it Suppose that the conditions of
Corollary 3.7 are satisfied. Then $G/_cH$ with the quotient topology
is discrete if and only if $H$ is open in $G$.}
\par The assertions of Corollaries 3.6-3.8 follow from Theorem 3.4.

\par {\bf Lemma 3.9.} {\it Suppose that $G$ is a topological loop,
$H$ is a closed subquasigroup in $G$ satisfying Condition $(3.2.1)$.
Then $cl_{G/_cH}\pi (V)\subset \pi (U/g_2)$ for each $g_2$ in $G$
and open subsets $U$ and $V$ in $G$ such that $g_2\in U$, $e\in V$
and $(V/V)g_2\subset U$.}
\par {\bf Proof.} We take any fixed $g$ in $G$ such that $\pi (g)\in
cl_{G/_cH}\pi (V)$, where the quotient map $\pi $ is as in Lemma
3.2. Notice that $gV$ is an open neighborhood of $g$ in $G$. Note
also that the quotient map $\pi : G\to G/_cH$ is open by Theorem
3.4. Therefore $\pi (gV)$ is an open neighborhood of $\pi (g)$ in
$G/_cH$, since $H(gV)=(Hg)V=\bigcup \{ (hg)V: h\in H \} $ is open in
$G$. Hence $\pi (gV)\cap \pi (V)\ne \emptyset $. This implies that
there exist $a$ and $b$ in $V$ with $\pi (ga)=\pi (b)$,
consequently, there exist $h_1$ and $h_2$ in $H$ with
$h_1(ga)=h_2b$.
\par Then we infer that $ga=h_1\setminus (h_2b)$. From $(3.2.1)$ it
follows that there exists $h_3\in H$ such that $h_1\setminus
(h_2b)=h_3b$, hence $g=(h_3b)/a$, consequently, there exists $h_4\in
H$ such that $g=h_4(b/a)\in H(b/a)$ by $(3.2.1)$. On the other hand,
$H(b/a)\subset H(V/V)\subset \pi (U/g_2)$, consequently, $\pi (g)\in
\pi (U/g_2)$. Thus $cl_{G/_cH}\pi (V)\subset \pi (U/g_2)$.

\par {\bf Theorem 3.10.} {\it Let $G$ be a topological $T_1$ loop, $H$ be a
closed subquasigroup of $G$ satisfying Condition $(3.2.1)$. Then the
quotient space $G/_cH$ is regular.}
\par {\bf Proof.} In view of Lemma 3.2 the quotient space $G/_cH$
and the quotient map $\pi : G\to G/_cH$ exist. Take any fixed
$g_2\in G$ and an open neighborhood $P$ of $\pi (g_2)$ in $G/_cH$.
From the continuity of the map $\pi $ it follows that there exists
an open neighborhood $U$ of $g_2$ in $G$ such that $\pi (U)\subset
P$. Since $G$ is the topological loop, then for each fixed $g_2\in
G$ there exists an open neighborhood $V$ of $e$ in $G$ such that
$V/V\subset U/g_2$. By virtue of Lemma 3.9 $cl_{G/_cH} \pi
(V)\subset \pi (U/g_2)$, that is $cl_{G/_cH} \pi (Vg_2)\subset \pi
(U)$. On the other hand, $\pi (Vg_2)$ is an open neighborhood of
$\pi (g_2)$ and $G/_cH$ is homogeneous and $T_1$ by Theorem 3.4.
This implies that $G/_cH$ is regular.

\par {\bf Definition 3.11.} If $X$ and $Y$ are topological spaces and
$f: X\to Y$ is a continuous map such that $f$ is closed and a
preimage $f^{-1}(y)$ is compact for each $y\in Y$, then $f$ is
called a perfect map.

\par {\bf Lemma 3.12.} {\it Assume that $G$ is a quasigroup, $A$, $B$,
$C$ are subsets in $G$. Then $(AB)\cap C=\emptyset $ if and only if
$A\cap (C/B)=\emptyset $ and if and only if $(A\setminus C)\cap
B=\emptyset $.}
\par {\bf Proof.} For each $a\in A$, $b\in B$, $c\in C$,
$(ab\ne c) \leftrightarrow (a\ne c/b) \leftrightarrow (a\setminus
c\ne b)$.

\par {\bf Theorem 3.13.} {\it Let $G$ be a topological $T_1$ loop, $F$ be a
compact subset in $G$, $P$ be a closed subset in $G$ and $F\cap
P=\emptyset $. Then an open neighborhood $V$ of $e$ in $G$ exists
such that $FV\cap P=\emptyset $ and $VF\cap P=\emptyset $.}
\par {\bf Proof.} The left translation $L_g: G\to G$ is the
homeomorphism of $G$ onto $G$. Therefore, for each $f\in F$ an open
neighborhood $V_f$ of $e$ in $G$ exists such that $fV_f\cap
P=\emptyset $. From the joint continuity of multiplication on $G$ it
follows that an open neighborhood $W_f$ of $e$ in $G$ exists with
$(fW_f)W_f\subset fV_f$. Hence $F\subset \bigcup_{f\in F} fW_f$ and
from the compactness of $F$ it follows that these open covering of
$F$ has a finite subcovering $F\subset \bigcup_{f\in \Lambda }
fW_f$, where $\Lambda \subset F$, $card (\Lambda )<\aleph _0$.
\par We put $S=\bigcap_{f\in \Lambda }W_f$, consequently, $S$ is
an open neighborhood of $e$ in $G$. For each $g\in F$ an element $f$
in $\Lambda $ exists with $g\in fW_f$, consequently, $gS\subset
(fW_f)W_f\subset fV_f\subset G-P$ and hence $FS\cap P=\emptyset $.
\par Symmetrically an open neighborhood $Q$ of $e$ in $G$ exists
such that $QP\cap F=\emptyset $. Taking $V=S\cap Q$, we get the
assertion of this theorem.

\par {\bf Theorem 3.14.} {\it Let $G$ be a topological $T_1$ loop, $F$ be a
compact subset in $G$, $P$ be a closed subset in $G$. Then the sets
$FP$ and $PF$ are closed in $G$.}
\par {\bf Proof.} In view of Lemma 3.12 for each $V\subset G$ the
condition $(Va)\cap (FP)=\emptyset $ is equivalent to $(V\setminus
(FP))\cap \{ a \} =\emptyset $. From the joint continuity of $Div_l$
and multiplication on $G$ it follows that for each $f\in F$ an open
neighborhood $V_f$ of $e$ in $G$ exists such that $(V_f\setminus
((V_ff)P))\cap \{ a \} =\emptyset $.
\par The subset $F$ in $G$ is compact, consequently, the open
covering $ \{ V_ff: f\in F \} $ of $F$ contains a finite subcovering
$ \{ V_ff: f\in \Lambda \} $, $F\subset \bigcup_{f\in \Lambda
}V_ff$, where $\Lambda $ is a finite subset in $F$. Therefore
$S=\bigcap_{f\in \Lambda }V_f$ is an open neighborhood of $e$ in
$G$. Hence $S\setminus (FP)\subset S\setminus ((\bigcup_{f\in
\Lambda }(V_ff))P)\subset \bigcup_{f\in F}V_f\setminus ((V_ff)P)$
and $(S\setminus (FP))\cap \{ a \} =\emptyset $, since $a\notin
\bigcup_{f\in F}V_f\setminus ((V_ff)P)$. Therefore $(Sa)\cap
(FP)=\emptyset $ by Lemma 3.12. Thus $FP$ is closed. Symmetrically
is proven that $PF$ is closed in $G$.

\par {\bf Theorem 3.15.} {\it Assume that $G$ is a topological $T_1$ loop
and $H$ is a compact subquasigroup in $G$ satisfying Condition
$(3.2.1)$. Then the quotient map $\pi : G\to G/_cH$ is perfect.}
\par {\bf Proof.} From Theorem 3.14 it follows that $HP$ is closed in
$G$ for each closed subset $P$ in $G$. On the other hand, we have
$HP=\pi ^{-1}(\pi (P))$ and $\pi (P)$ is closed in $G/_cH$,
consequently, the quotient map $\pi $ is closed. If $y\in G/_cH$ and
$g\in G$ with $\pi (g)=y$, then $\pi ^{-1}(y)=Hg$ is compact in $G$,
since $H$ is compact and $R_g: G\to G$ is the homeomorphism of $G$.
Hence the fibers of $\pi $ are compact. Thus the quotient map $\pi $
is perfect.

\par {\bf Corollary 3.16.} {\it Suppose that the conditions of Theorem
3.15 are satisfied and $G/_cH$ is compact. Then $G$ is compact.}
\par {\bf Proof.} According to Theorem 3.15 the quotient map $\pi :
G\to G/_cH$ is perfect. Since $G/_cH$ is compact, then $G$ is
compact by Theorem 3.7.2 in \cite{eng}.

\par {\bf Corollary 3.17.} {\it Let the conditions of Corollary
1 in \cite{ludstwpmgax19} be satisfied and let ${\cal C}_1$ and $A$
and $A{\cal C}_1$ be closed in $G$, then $D/_c(A{\cal C}_1)$ and
$(A{\cal C}_1)/_cA$ are homogeneous $T_1\cap T_3$ spaces and the
quotient maps $\pi ^D_{A{\cal C}_1}: D\to D/_c(A{\cal C}_1)$ and
$\pi ^{A{\cal C}_1}_A: A{\cal C}_1\to (A{\cal C}_1)/_cA$ are open.}
\par {\bf Proof.} This follows from Theorems 3.4, 3.10 and Example 3.3.

\par {\bf Corollary 3.18.} {\it Let the conditions of Theorem 2.3 be
satisfied and let $(G, {\cal T}_{G,G_1,G_2})$ be compact and
$(G,{\cal T}_G)$ be $T_1$ as the topological quasigroup. Then ${\cal
T}_{G,G_1,G_2}={\cal T}_G$.}
\par {\bf Proof.} Since $(G,{\cal T}_G)$ is the $T_1$ topological
quasigroup, then it is regular. In view of Corollary 3.1.14 in
\cite{eng}  ${\cal T}_{G,G_1,G_2}={\cal T}_G$.

\par {\bf Corollary 3.19.} {\it Assume that the conditions of
Corollary 2.4 are satisfied and $(D, {\cal T}_{D,A{\cal C}_1,A})$ is
compact. Then ${\cal T}_D={\cal T}_{D,A{\cal C}_1,A}$; moreover, $A$
and $A{\cal C}_1$ are compact relative to the topologies ${\cal
T}_D\cap A$ and ${\cal T}_D\cap (A{\cal C}_1)$ respectively
inherited from $G$.}
\par {\bf Proof.} Corollary 3.18 implies that ${\cal T}_D={\cal
T}_{D,A{\cal C}_1,A}$. The maps $\tau ^D_A$ and $\tau ^D_{A{\cal
C}_1}$ are continuous by Corollary 2.4. On the other hand, $A=(\tau
^D_A)^{-1}(e)$ and $A{\cal C}_1=(\tau ^D_{A{\cal C}_1})^{-1}(e)$,
consequently, $A$ and $A{\cal C}_1$ are closed in $D$, hence $A$ and
$A{\cal C}_1$ are compact by Theorem 3.1.2 in \cite{eng}.

\par {\bf Proposition 3.20.} {\it Assume that the conditions of
Remark 1 in \cite{ludstwpmgax19} are satisfied, $G=A\star ^{\phi ,
\eta , \kappa , \xi }B$ is a smashed twisted product of metagroups
$A$ and $B$ with smashing factors $\phi $, $\eta $, $\kappa $, $\xi
$. Then embeddings $\theta ^G_A: A\hookrightarrow G$ and $\theta
^G_B: B\hookrightarrow G$ exist, and $\theta ^G_B(B)$ in $G$ is
invariant. Moreover, a transversal set $V_{G,B}$ exists such that
$V_{G,B}=\theta ^G_A(A)$.}
\par {\bf Proof.} We shortly denote $\theta ^G_A$ as $\theta _A$,
because $G$ is specified, and we put $\theta _A(a)=(a,e)$ with
$e=e_B$ for each $a\in A$; $\theta _B(b)=(e,b)$ with $e=e_A$ for
each $b\in B$. From Formula $(37)$ in \cite{ludstwpmgax19} it
follows that $(e,b)(a,e)=(a,b\xi ((e,b), (a,e)))$ for each $a\in A$
and $b\in B$. Therefore, for each $g=(a_1,b_1)$ in $G$ there exist
unique $a\in A$ and $b\in B$ such that \par $(3.20.1)$
$(e,b)(a,e)=g$ with $a=a_1$ and $b=b_1/\xi ((e,b_1), (a_1, e))$, \\
since $\xi ((e,b), (a,e))=\xi ((e,b_1), (a,e))$ by $(35)$ in
\cite{ludstwpmgax19}. Certainly the maps $a_1\mapsto a$ and
$(a_1,b_1)\mapsto b$ provided by $(3.20.1)$ are single-valued.
\par For each $g_1=(e,b_1)\in G$, $g_2=(a_2,b_2)\in G$,
$g_3=(a_3,b_3)\in G$ we deduce that \par
$I_1=(g_1g_2)g_3=(a_2a_3,b_3^{a_2}(b_2b_1)\xi ((e,b_1),(a_2,b_2))
\xi ((a_2,b_2b_1), (a_3,b_3)))$ and
\par $I_2=g_1(g_2g_3)=(a_2a_3, (b_3^{a_2}b_2) \xi ((a_2,b_2),
(a_3,b_3)) b_1 \xi ((e,b_1), (a_2a_3,b_3^{a_2}b_2)))$ \\ by
Conditions $(31)$, $(32)$, $(34)$ in \cite{ludstwpmgax19}. Hence
$I_1=tI_2$ with $t=t(g_1,g_2,g_3)\in \theta _B({\cal C})$,
consequently, $\theta _B(B)$ satisfies Condition $(3.2.1)$, since
${\cal C}\hookrightarrow {\cal C}(B)$ by Remark 1 in
\cite{ludstwpmgax19}.
\par In view of Lemma 3.2 and Formula $(3.20.1)$ the
transversal set $V_{G,B}=\theta _A(A)$ and the maps $\psi =\psi
^G_B: G\to \theta _B(B)$ and $\tau =\tau ^G_B: G\to \theta _A(A)$
exist such that \par $(3.20.2)$ $g=g^{\psi }g^{\tau }$ with
\par $(3.20.3)$ $g^{\psi }=(e,b)$ and $g^{\tau }=(a,e)$ for each
$g=(a_1,b_1)\in G$, where \par $(3.20.4)$ $a=a_1$ and $b=b_1/\xi
((e,b_1), (a_1,e))$.

\par It remains to prove that $\theta _B(B)$ is invariant in $G$.
For this it is sufficient to prove that
\par $(3.20.5)$ $g_1\theta _B(B)=\theta _B(B)g_1$ and
\par $(3.20.6)$ $(g_1\theta _B(B))g_2=g_1(\theta _B(B)g_2)$
for each $g_1$ and $g_2$ in $G$, \\ since Properties $(3.20.5)$ and
$(3.20.6)$ imply that $(g_1g_2)\theta _B(B)=g_1(g_2\theta _B(B))$
for each $g_1$ and $g_2$ in $G$.
\par For each $g_1=(a_1,b_1)$ in $G$ and $b_2\in B$, $b_3\in B$ we
get that \par $(a_1,b_1)(e,b_2)=(a_1,b_2^{a_1}b_1\xi ((a_1,b_1),
(e,b_2)))$ and \par $(e,b_3)(a_1,b_1)=(a_1, b_1b_3 \xi ((e,b_3),
(a_1,b_1)))$ \\ according to $(37)$ in \cite{ludstwpmgax19}. The
following equation
\par $b_2^{a_1}b_1\xi ((a_1,b_1), (e,b_2)) =b_1b_3 \xi ((e,b_3),
(a_1,b_1))$ has a unique solution
\par $b_3=[b_1\setminus (b_2^{a_1}b_1)]\xi ((a_1,b_1), (e,b_2))/\xi ((e,[b_1\setminus (b_2^{a_1}b_1)]),
(a_1,b_1))$ for given $g_1=(a_1,b_1)$ and $g_2=(e,b_2)$, since $\xi
$ satisfies Condition $(35)$ in \cite{ludstwpmgax19}. From $\xi
(g_1,g_3)\in {\cal C}$ for each $g_1$ and $g_3$ in $G$ and ${\cal
C}\hookrightarrow {\cal C}(B)\subset B$ it follows that $b_3\in B$.
Thus $G$ satisfies Condition $(3.20.5)$.
\par Then we consider $I_1=(g_1(e,b_2))g_3$ and
$\tilde{I}_2=g_1((e,\tilde{b}_2)g_3)$ for any $g_1$ and $g_3$ in
$G$, $b_2$ and $\tilde{b}_2$ in $B$. Then we infer that
\par $I_1=(a_1a_3, b_3^{a_1}(b_2^{a_1}b_1) \xi ((a_1,b_1), (e,b_2))
\xi ((a_1,b_2^{a_1}b_1), (a_3,b_3)))$ and
\par $\tilde{I}_2=(a_1a_3, (b_3^{a_1}\tilde{b}_2^{a_1}) \kappa (a_1,
b_3, \tilde{b}_2) \xi ((e, \tilde{b}_2), (a_3,b_3)) b_1 \xi
((a_1,b_1), (a_3, b_3 \tilde{b}_2)))$ \\ by $(33)$  and $(37)$
\cite{ludstwpmgax19}. The following equation $I_1=\tilde{I}_2$ is
satisfied if and only if
\par $b_2^{a_1}b_1\gamma =\tilde{b}_2^{a_1}b_1\tilde{\alpha }
p(b_3^{a_1},\tilde{b}_2^{a_1},b_1\tilde{\alpha })$ with \par $\gamma
= \xi ((a_1,b_1), (e,b_2)) \xi ((a_1,b_2^{a_1}b_1), (a_3, b_3))$;
\par $\tilde{\alpha }=\kappa (a_1,b_3,\tilde{b}_2)
\xi ((e,\tilde{b}_2),(a_3,b_3)) \xi ((a_1,b_1), (a_3,
b_3\tilde{b}_2))$ \\ by $(1.1.1)$ and $(1.1.9)$. Using $(34)$,
$(35)$ and Lemma 2 in \cite{ludstwpmgax19}  we deduce that there
exists a unique solution
\par $(3.20.7)$ $\tilde{b}_2^{a_1}=(b_2^{a_1}b_1\gamma )/(b_1\delta )$
for the given $g_1=(a_1,b_1)$ and $g_3=(a_3,b_3)$ in $G$, $b_2\in B$
with \par $\delta = \alpha p(b_3^{a_1}, b_2^{a_1}, b_1\alpha )$ and
\par $\alpha =\kappa (a_1,b_3,b_2) \xi
((e,b_2),(a_3,b_3)) \xi ((a_1,b_1), (a_3, b_3b_2))$.
\par Since $\gamma \in {\cal C}$ and $\delta \in {\cal C}$, ${\cal C}\hookrightarrow
{\cal C}(B)$, then $\tilde{b}_2^{a_1}\in B$. From $(31)$ and $(33)$
\cite{ludstwpmgax19} it follows that
\par $(3.20.8)$ $\tilde{b}_2=(\tilde{b}_2^{a_1})^{e/a_1}/ \eta (e/a_1,
a_1, (\tilde{b}_2^{a_1})^{e/a_1})$. \\ Hence $(3.20.7)$ and
$(3.20.8)$ imply that $\tilde{b}_2\in B$, consequently, $G$
satisfies Condition $(3.20.6)$. Thus $\theta _B(B)$ is the invariant
submetagroup in $G$.

\par {\bf Corollary 3.21.} {\it If the conditions of Remark 1
in \cite{ludstwpmgax19} are satisfied, $A$ and $B$ are topological
$T_1$ metagroups, the topology on $G$ is induced by the Tychonoff
product topology on $A\times B$ and the smashing factors $\phi $,
$\eta $, $\kappa $, $\xi $ are (jointly) continuous, then the maps
$\psi : G\to \theta _B(B)$ and $\tau : G\to V_{G,B}=\theta _A(A)$
are continuous relative to the topology ${\cal T}_G$ on the
topological metagroup $G=A\star ^{\phi , \eta , \kappa , \xi }B$.}
\par {\bf Proof.} This follows from Formulas $(3.20.2)$-$(3.20.4)$ and
the (joint) continuity of the smashing factors $\phi $, $\eta $,
$\kappa $, $\xi $ and hence of $Div_r$ and $t_G$ on $(G, {\cal
T}_G)$, where the topology ${\cal T}_G$ on $G$ is induced by the
Tychonoff product topology on $A\times B$.

\par {\bf Corollary 3.22.} {\it Let for pairs of metagroups
$A_j$, $B_j$ the conditions of Remark 1 in \cite{ludstwpmgax19} be
satisfied for each $j\in \{ 1, 2 \} $, where $B_1=B_2$ such that
${\cal C}_m(A_j)\subset {\cal C}\hookrightarrow B_j\hookrightarrow
{\cal C}(A_j)$ for each $j\in \{ 1, 2 \} $. Let $\phi
_2(a)b=id(b)=b$ for each $a\in {\cal C}$ and $b\in B_2$, and $\xi
_2((a,e),(e,b))=\xi _2((e,b),(a,e))$ for each $a\in A_2$ and $b\in
B_2$. Let $A'=A_1\star ^{\phi _1, \eta _1, \kappa _1, \xi _1}B_1$
and $B=A_2\star ^{\phi _2, \eta _2, \kappa _2, \xi _2}B_2$ and let
$\phi _3$, $\eta _3$, $\kappa _3$, $\xi _3$ for the pair $(A',B)$
with ${{\cal C}'}_1=\theta ^B_{B_2}(B_2)$ satisfy the conditions of
Remark 1 in \cite{ludstwpmgax19} (with ${{\cal C}'}_1$ instead of
${\cal C}$) and let $D=A'\star ^{\phi _3, \eta _3, \kappa _3, \xi
_3}B$, where $\theta ^B_{B_2}: B_2\hookrightarrow B$ is the
embedding provided by Proposition 3.20. Then there are embeddings
$\theta _{A_j}: A_j\hookrightarrow D$, $\theta _{B_j}:
B_j\hookrightarrow D$ for each $j\in \{ 1, 2 \} $, $\theta _B:
B\hookrightarrow D$ such that $D$ with $A=\theta _2(A_2)$ and ${\cal
C}_1=\theta _{B_2}(B_2)$ satisfy Condition $(28)$ in
\cite{ludstwpmgax19} and $xA=Ax$ for each $x\in {\cal C}_1$.}
\par {\bf Proof.} By virtue of Theorem 4 in \cite{ludstwpmgax19} $B$,
$A'$ and $D$ are metagroups and there are embeddings $\theta _{A_j}:
A_j\hookrightarrow D$, $\theta _{B_j}: B_j\hookrightarrow D$ for
each $j\in \{ 1, 2 \} $, $\theta _B: B\hookrightarrow D$ such that
$\theta _B(B)=\theta _{A_2}(A_2) \theta _{B_2}(B_2)$, since
$B_2\hookrightarrow {\cal C}(A_2)$.
\par For each $(e,b)\in B$, $(a,e)\in B$ and $(a_2,e)\in B$ with
$b\in B_2$, $a\in A_2$, $a_2\in A_2$ we deduce that
\par $(a,e) (e,b)=(a,b^a \xi ((a,e), (e,b)))$ and
\par $(e,b) (a_2,e) = (a_2, b \xi ((e,b), (a_2,e)))$.
\\ Therefore $(a,e) (e,b)=(e,b) (a_2,e)$ if and only if
$a=a_2$ and $b^a \xi ((a,e), (e,b))=b \xi ((e,b), (a_2,e))$. From
$(31)$ and $(35)$ in \cite{ludstwpmgax19} and the conditions of this
corollary it follows that $xA=Ax$ for each $x\in {\cal C}_1$, since
$\phi_2(a)b=b^a=b$ for each $a\in A_2$ and $b\in B_2$. In view of
Proposition 3.20 the subgroup ${\cal C}_1$ is invariant in $\theta
_B(B)$ and $\theta _{A'}(A')$.
\par Certainly, $\theta
_{B_1}(B_1)$ and $\theta _{B_2}(B_2)$ are isomorphic subgroups in
$D$, since $B_1=B_2$. Hence each $d\in D$ can be presented in the
following form: $d=a_1a_2b$ with $a_1\in \theta _{A_1}(A_1)$,
$a_2\in \theta _{A_2}(A_2)$ and $b\in \theta _{B_1}(B_1)$.  From
$a_1\theta _{B_2}(B_2)=\theta _{B_2}(B_2)a_1$ and $a_2\theta
_{B_2}(B_2)=\theta _{B_2}(B_2)a_2$ it follows that $d\theta
_{B_2}(B_2)=\theta _{B_2}(B_2)d$ for each $d\in D$. On the other
hand, ${\cal C}_m(D)\subset {\cal C}_1$, since ${\cal
C}_m(A_j)\subset {\cal C}\hookrightarrow B_j\hookrightarrow {\cal
C}(A_j)$ for each $j\in \{ 1, 2 \} $. Consequently, the subgroup
$\theta _{B_2}(B_2)$ is invariant in $D$.

\par {\bf Corollary 3.23.} {\it Assume that the conditions of
Corollary 3.22 are satisfied, $A_j$, $B_j$ are $T_1$ topological
metagroups for each $j\in \{ 1, 2 \} $ and $\phi _i$, $\eta _i$,
$\kappa _i$, $\xi _i$ are jointly continuous for each $i\in \{ 1, 2,
3 \} $. Then $D$, $A$, ${\cal C}_1$ provided by Corollaries 3.21 and
3.22 are $T_1\cap T_3$ topological metagroups and satisfy the
conditions of Theorem 6 in \cite{ludstwpmgax19} and ${\cal C}_1$ is
closed in $D$.}
\par {\bf Proof.} This follows from Theorem 4 in
\cite{ludstwpmgax19} and Corollaries 3.21 and 3.22.

\par {\bf Definition 3.24.} Let $\Lambda $ be a directed set, $G_j$ be
a topological metagroup (or loop), $\pi ^j_i: G_j\to G_i$ be a
continuous homomorphism for each $i\le j$ in $\Lambda $ such that
$\pi ^j_i\circ \pi ^k_j=\pi ^k_i$ for each $i\le j\le k$ in $\Lambda
$, and $\pi ^i_i=id_{G_i}$ for each $i\in \Lambda $, where
$id_{G_i}(g_i)=g_i$ for each $g_i\in G_i$. Then $S= \{ G_j, \pi
^j_i, \Lambda \} $ is called an inverse continuous homomorphism
system of topological metagroups (or loops respectively). If a
topological metagroup $G$ is a limit of $S$, $G=\overleftarrow{\lim
S}$, then it is said that $G$ is decomposed into $S$.

\par {\bf Theorem 3.25.} {\it There exists an infinite family ${\cal
F}$, each $G\in {\cal F}$ is a topological $T_1$ metagroup such that
$G$ is compact, locally connected and can not be decomposed into the
inverse continuous homomorphism system $S_G=\{ G_j, \pi ^j_i,
\Lambda \} $ of topological metagroups $G_j$ with $dim (G_j)<\infty
$ for each $j\in \Lambda $.}
\par {\bf Proof.} We take any locally connected $T_1$ compact metagroups $A$,
$B$ and their invariant closed subgroup ${\cal C}$ with positive
covering dimensions $dim (A)>0$, $dim (B)>0$, $dim (A/_c{\cal
C})>0$, $dim (B/_c{\cal C})>0$ such that the conditions of Theorem 6
in \cite{ludstwpmgax19} are satisfied. Evidently, such triples $(A,
B, {\cal C})$ exist and their family is infinite. Indeed, in
particular, they may be direct products $A=K_1\times {\cal C}$,
$B=P_1\times {\cal C}$ or semidirect products $A=K_1\times ^s{\cal
C}$, $B=P_1\times ^s{\cal C}$ with topological $T_1$ metagroups
$K_1$, $P_1$ and a topological $T_1$ group ${\cal C}$; or
particularly $A$, $B$ may be topological $T_1$ groups.
\par Therefore, $\theta _B(B)$ is invariant in the smashed
twisted product $G=A\star ^{\phi , \eta , \kappa , \xi }B$ such that
$G$ is a topological $T_1$ metagroup and a transversal set exists
$V_{G,B}=\theta _A(A)$ by Corollary 1 in \cite{ludstwpmgax19} and
Proposition 3.20. By virtue of Corollary 3.21 the maps $\psi : G\to
\theta _B(B)$ and $\tau : G\to V_{G,B}$ are continuous. For compact
$A$ and $B$ the metagroup $G$ is compact by the Tychonoff theorem
3.2.4 in \cite{eng}.
\par This implies that there are triples $(A_1,B_1,{\cal C})$ and $(A_2,B_2,{\cal C})$
satisfying the conditions of Corollary 3.23 with locally connected
$T_1$ compact metagroups $A_1$, $A_2$, $B_1=B_2$ and their invariant
closed subgroup ${\cal C}$ with positive covering dimensions $dim
(A_1)>0$, $dim (A_2)>0$, $dim (B_1)>0$, $dim (A_1/_c{\cal C})>0$,
$dim (A_2/_c{\cal C})>0$, $dim (B_1/_c{\cal C})>0$. Then $D$, $A$,
${\cal C}_1$ provided by Corollaries 3.21 and 3.22 satisfy the
conditions of Corollary 1 in \cite{ludstwpmgax19} or Proposition
3.20 such that $A=\theta _{A_2}(A_2)$. By virtue of Theorem 6 in
\cite{ludstwpmgax19}, Theorem 1 in \cite{ludaxsftm20}, Theorem 3.15
and Corollary 3.23 $D$, $A$, $B$ are locally connected $T_1$ compact
metagroups with a closed invariant subgroup ${\cal C}_1=\theta
_{B_2}(B_2)$ and $dim (D)>0$, $dim (A)>0$, $dim (B)>0$. Moreover,
$V_{D,A{\cal C}_1}= \theta _{A'}(A')$, with $A'=A_1\star ^{\phi _1,
\eta _1, \kappa _1, \xi _1} B_1$, and there is a bijection from
$V_{A{\cal C}_1,A}=V_{{\cal C}_1,{\cal C}_{1,A}}$ onto $(A{\cal
C}_1)/_cA$ by Remark 2.6 and $V_{{\cal C}_1,{\cal
C}_{1,A}}V_{D,A{\cal C}_1}=V_{D,A}$ by Formula $(2.6.3)$. Therefore,
$V_{{\cal C}_1,{\cal C}_{1,A}}$ and $V_{D,A{\cal C}_1}$ can be
chosen compact, consequently, $V_{D,A}$ is compact by the Tychonoff
theorem 3.2.4 in \cite{eng}. Corollaries 3.17 and 3.21 imply that
the maps $\psi ^D_A: D\to A$ and $\tau ^D_A: D\to V_{D,A}$ are
continuous relative to the topology ${\cal T}_D$.
\par In view of Theorems 2.10 and 2.11 there exists a $T_1\cap T_3$ compact metagroup
$C_0=D\Delta _A^{\phi , \eta , \kappa , \xi }F_0$, where $F_0$ is
closed in $(C(V,B),{\cal T}_{\cal W})$, where $F_0\subset
C(V,B)\subset F=B^V$, $V=V_{D,A}$. Hence $dim (C_0)>0$ and $C_0$ is
locally connected. In view of Theorem 3.1.9 \cite{eng} $D$, $A$,
$B$, ${\cal C}_1$, $C_0$ are $T_1\cap T_4$ topological spaces. By
the construction above $card (V)\ge \aleph _0$. \par On the other
hand, a family $Hom_{c,{\cal C}}=Hom_{c,{\cal C}}((A\times
F_0)\times (A\times F_0),{\cal C}_1)$ of all continuous
homomorphisms from $(A\times F_0)\times (A\times F_0)$ into ${\cal
C}_1$ satisfying $(35)$ in \cite{ludstwpmgax19} is a proper closed
subset in a family $C_{\cal C}=C_{\cal C}((A\times F_0)\times
(A\times F_0),{\cal C}_1)$ of all continuous maps from $(A\times
F_0)\times (A\times F_0)$ into ${\cal C}_1$ satisfying $(35)$ in
\cite{ludstwpmgax19}. Since $dim (B)>0$ and $dim (V)>0$, then $dim
(C(V,B))=\infty $, where $C(V,B)$ is in the ${\cal T}_{\cal W}$
topology.
\par We choose $F_0$ with $dim (F_0)=\infty $ and the map $\xi \in C_{\cal C}-
Hom_{c,{\cal C}}$ with values in ${\cal C}_1$ such that $\xi
((d_1^{\psi },f_1), (d^{\psi },f))(v)$ depends nontrivially on
infinite number of coordinates $v\in V$ for an infinite family of
$(f_1,f)\in F_0\times F_0$ for each $d_1^{\psi }\ne e$ and $d^{\psi
}\ne e$, where $d$ and $d_1$ belong to $D$, $f\in F_0$, $f_1\in
F_0$, since $dim ({\cal C}_1)>0$ and $dim (A_2)>0$. This implies
that there exists the topological metagroup $G=C_0$ with $dim
(G)=\infty $ which can not be decomposed into the inverse continuous
homomorphism system $S_{G}=\{ G_j, \pi ^j_i, \Lambda \} $ of
topological metagroups $G_j$ with $dim (G_j)<\infty $ for each $j\in
\Lambda $. From the proof above it follows that the family of such
topological metagroups is infinite.
\par {\bf Remark 3.26.} If instead of Theorem 6 use Theorem 5 in
\cite{ludstwpmgax19}, then Theorem 3.25 will be for topological
loops $G$.
\par {\bf 3.27. Conclusion.}  The obtained results in this paper can be used for
subsequent investigations of topological metagroups, quasigroups,
loops, topological algebras, generalized $C^*$-algebras, operator
algebras over octonions and the Cayley-Dickson algebras and
noncommutative geometry associated with them \cite{bruckb}.
Certainly, metagroups can be realized using deterministic functions,
which are important in algorithm theory and informatics. Besides
these applications of metagroups and those outlined in the
introduction, it is interesting to indicate possible applications in
mathematical coding theory and related reliability of software and
hardware systems \cite{blautrctb,petushrtj820,srwseabm14}, because
frequently codes are based on topological binary systems. \\
\par {\bf Declarations.}
\par {\bf Funding.} This research has not received any funding.
\par {\bf Conflicts of interests/Competing interests.} The authors
declare to have no any conflict of interests.

\end{document}